\documentclass[12pt]{article}

 \usepackage{graphics}
 \usepackage{graphicx}
 \usepackage{epsfig}
 \usepackage[all]{xy}

 \usepackage{amssymb}
 \usepackage{amsthm,amsmath}

 \usepackage{latexsym, epsfig, psfrag,eepic,colordvi,bm}
 \usepackage{graphics,color}
 \usepackage{amsmath,amsfonts,amssymb,amscd}
 \usepackage[all]{xy}
 \usepackage{epstopdf}

 \usepackage{amsthm,amsmath,amssymb}

 \usepackage{graphicx}
 \usepackage{subcaption}
\usepackage{amsthm,amsmath,amssymb}
\usepackage{graphicx,epsfig}

\usepackage[utf8]{inputenc}

\usepackage{algorithm}
\usepackage{algorithmicx}
\usepackage[noend]{algpseudocode}

\makeatletter
\def\BState{\State\hskip-\ALG@thistlm}
\makeatother

\usepackage[table]{xcolor}

\newcolumntype{s}{>{\columncolor[HTML]{AAACED}} p{3cm}}

\theoremstyle{plain}
\newtheorem{theorem}{Theorem}[section]

\theoremstyle{definition}

\theoremstyle{remark}

\DeclareMathOperator{\Deg}{deg}

\date{}

\title{D-chain tomography of networks: a new structure spectrum and an application to the SIR process}

\author{Ricky X. F. Chen, Christian M. Reidys, Andrei C. Bura}

\begin{document}

\maketitle
\noindent{\bf Abstract.}
The analysis of the dynamics on complex networks is closely connected to structural features of
the networks. Features like, for instance, graph-cores and node degrees have been studied ubiquitously.
Here we introduce the D-spectrum of a network, a novel new
framework that is based on a collection of nested chains of subgraphs within the network. Graph-cores and node degrees are merely from two particular such chains of the D-spectrum. Each chain
gives rise to a ranking of nodes and, for a fixed node, the collection of these ranks provides us with
the D-spectrum of the node. Besides a node deletion algorithm, we discover a connection between the D-spectrum of a network and some fixed points of certain graph dynamical systems (MC systems) on the network. Using the D-spectrum we
identify nodes of similar spreading power in the susceptible-infectious-recovered (SIR) model on
a collection of real world networks as a quick application. We then discuss our results and conclude that D-spectra represent
a meaningful augmentation of graph-cores and node degrees.



\section{Introduction}
Structural properties of complex networks are of central importance for understanding the formation
principles of said networks and dynamics associated to them. Various network features have been
studied, such as node-degree~\cite{scale,scale2}, path distance~\cite{path,path2}, $k$-core
decomposition~\cite{kcore1,kcore2,kcore3}, motif identification~\cite{motif} and community
identification~\cite{commu1,commu2}.
Spreading dynamics on networks, such as, for instance, information diffusion, knowledge
dissemination, disease spreading, etc.~\cite{spread1,spread2,spread3,spread4,spread5,spread6,spread7}
is ubiquitous and has been studied extensively.
In the analysis particular focus has been put on the identification of nodes, that are the most
effective spreaders~\cite{import1,import2,import3,import5,import6}. Their localization
is of key relevance for designing strategies to decelerate or stop the spread, for instance in
infectious disease outbreaks, or accelerate the process as is in the case of knowledge dissemination. 

At first glance, the most connected nodes (hubs) seem to be natural candidates for being ``good'' spreaders.
However, Kitsak et al.~\cite{np-core} argue that the ``location'' of a node is more important than its
degree, where said location is characterized by its core number~\cite{kcore-b}.
These two perspectives differ significantly in that degree is a local feature while graph-cores are
(potentially) extended subgraphs. Recently, h-index families were proposed as a measure, and it was
shown that h-index outperforms both, degree as well as core-based measures in several
cases~\cite{nc-core}. In addition, discussion of an integration of node degree, h-index as well as
core number was presented there; a line of thought that can also be found implicitly in an earlier
work of Montresor, Pellegrini and Miorandi~\cite{mont}.

In this paper we present a completely new framework of charaterizing network structure
by introducing D-spectrum of a network. As a consequence, we obtain D-spectra of nodes which integrates node degrees and core numbers as
 endpoints of a sequence
which represent a transition from local to global information.
We present two methods to derive the D-spectra of nodes for a given network. The first
is a node deletion algorithm while the second is facilitated by computing fixed points of specific graph dynamical systems
on the network which we call
MC systems. 

As a quick application aiming at showing the potential practical value of the new framework, we then employ the D-spectra of nodes in order to characterize node similarity in the context of the
susceptible-infectious-recovered (SIR) model. We evaluate our approach based on the data from five distinct
networks and observe that  D-spectra represent
a meaningful augmentation of graph-cores and node degrees.  In the following we shall use the notions of graph and network interchangeably, as well as
those of node and vertex.

\section{D-spectrum}
In this section, we present the framework of D-chains of networks, the theory on MC systems, and their connection.
\subsection{D-chains of networks}
Suppose $G$ is a graph (without loops and multiple edges for simplicity). We write $H\leq G$ if $H$ is a
subgraph of $G$, and write $H<G$ if $H\leq G$ but
$H\neq G$. 

Let $L\colon G_0\geq  G_1 \geq  G_2 \geq  \cdots \geq G_k$ be a chain of nonempty subgraphs of
$G$, where $G_0=G$, and $G_i$ is a vertex-induced subgraph of $G_{i-1}$ for $1\leq i \leq k$.  The chain
$L$ is called a D-chain of order $t$ if for any $0\leq i \leq k$, every vertex $v\in G_i$ has at least
$i$ neighbors in $G_j$ where $j=\max \{0, i+t\}$. We call the number $k$ the length of the chain $L$
and denote $|L|=k$.
Clearly, each graph $G$ has a D-chain of order $t$ for any non-positive integer $t$, since $G_0=G$ is
a D-chain of order $t$ of length $0$. 

\vskip 12pt

\begin{equation*}
\xymatrix{
	G_0  & \cdots & G_{i-1+t} & G_{i+t} & \cdots & G_{i-1}\ar@/_2pc/[lll] & G_i\ar@/_2pc/[lll] & \cdots
}
\end{equation*}

A D-chain of order $t$, $L\colon G_0 \geq  G_1 \geq G_2 \geq
\cdots \geq G_k$, is called maximal if
(i)  there does not exist a D-chain $L'$ with $|L'|>k$; and
(ii) there does not exist a D-chain $L': G'_0\geq G'_1 \geq G'_2 \geq \cdots \geq G'_k$, where for
     some $ 1\leq i \leq k$, $G_i < G'_i$.  
For any $G$ and non-positive integer $t$, there exists a unique maximal D-chain of order $t$, because
the union of two D-chains of order $t$ of maximum length is again a D-chain of order $t$ of the same
length (see the Supplementary Information, Proposition C.1). 
Maximal D-chains are related as follows, where $G\rightarrow H$ denotes $H$ being a subgraph of $G$:
\begin{equation*}
\xymatrix{
	& \vdots \ar[d] & \vdots \ar[d] & \vdots \ar[d] & \\
	\cdots\ar[r] & G_{i-1}^{t-1} \ar[r]\ar[d] & G_{i}^{t-1} \ar[r]\ar[d] & G_{i+1}^{t-1} \ar[r]\ar[d] & \cdots\\
		\cdots \ar[r]& G_{i-1}^{t} \ar[r]\ar[d] & G_{i}^{t} \ar[r] \ar[d]& G_{i+1}^{t} \ar[r]\ar[d] & \cdots\\
			\cdots \ar[r]& G_{i-1}^{t+1} \ar[r] \ar[d]& G_{i}^{t+1} \ar[r] \ar[d]& G_{i+1}^{t+1} \ar[r] \ar[d]& \cdots\\
			& \vdots & \vdots & \vdots & 	
	}
\end{equation*}

Let $L\colon G_0 \geq  G_1 \geq G_2 \geq \cdots \geq G_k$ be the maximal D-chain of order $t$ of $G$.
We next introduce the rank of a node, $C_t$, by setting $C_t(v)=i$ if and only if $v$ is contained in
$G_i$ but not contained in $G_{i+1}$.
Let $\Delta(G)=\max_{v\in V(G)}\{\Deg(v)\}$, where $V(G)$ denotes the vertex set of $G$ and $\Deg(v)$
denotes the degree of $v$. We call the vector $(C_0(v), C_{-1}(v),\ldots, C_{-\Delta(G)}(v))$ the
D-spectrum of the vertex $v$. The collection of all maximal D-chains of $G$ is called the D-spectrum
of $G$. The
D-spectrum includes the chains for $t=0$ and $t=-\Delta(G)$, which produce the nested sequences of
$k$-cores and vertex degrees, respectively (see the Supplementary Information).

In Figure~$1$ we display specific maximal $D$-chains for a small network, and their embeddings
in the latter. We also display the D-spectra of all nodes of the network. The induced rank of a vertex is
represented by its color. 

	\begin{figure}[!htb]\label{F:jkj}
	\setlength{\tabcolsep}{0pt}
			\centering
	\begin{tabular}{@{\hspace{0em}}p{0.6\textwidth}@{\hspace{-1em}}p{0.3\textwidth}}
		\centering
		\begin{subfigure}{0.6\textwidth}
			\caption{}
			\centering
			\includegraphics[width=0.7\linewidth]{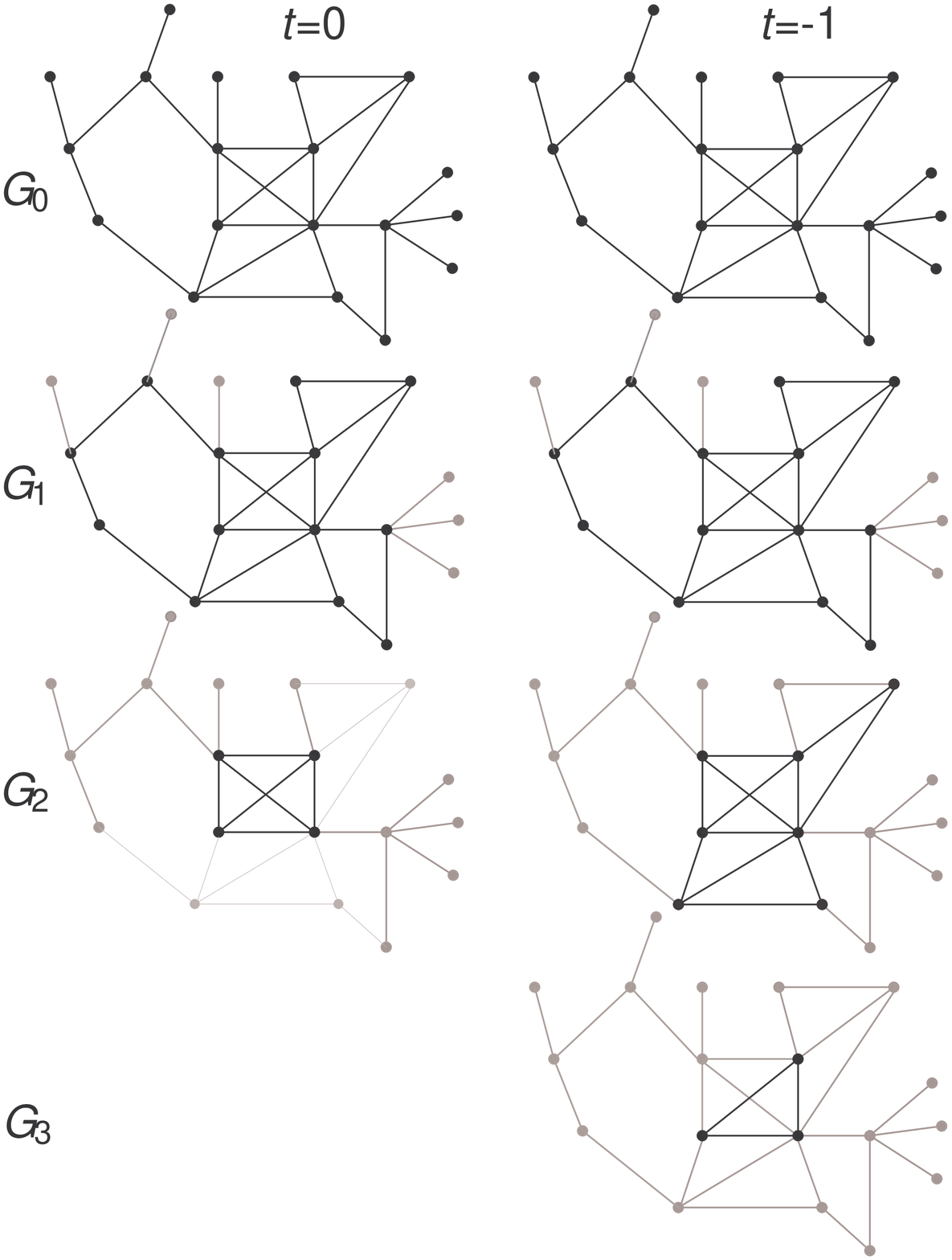}
		\end{subfigure}
		&	\begin{subfigure}{0.3\textwidth}
			\caption{}
			\centering
			\includegraphics[width=0.9\linewidth]{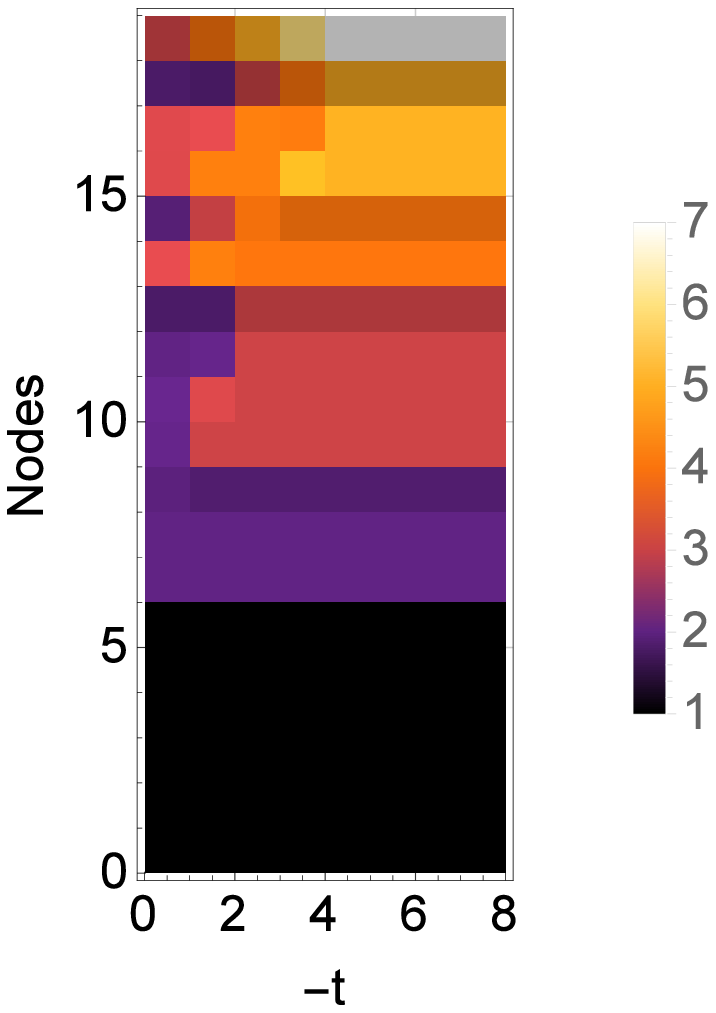}
		\end{subfigure}
	\end{tabular}
	\caption{(a) the maximal D-chains of order $t=0$ and $t=-1$ are highlighted in black. (b) the D-spectra
          of all vertices, where the induced ranks are indicated by colors.}
\end{figure}

\subsection{The D-spectrum via a deletion algorithm}

For fixed $k>0$ and $t<0$, suppose $i-mt=k$ for some $m\geq 0$ and $1\leq i \leq -t$. Given a graph $G$, we shall show that the following algorithm (Algorithm~$1$) will give us the maximal D-chain of order $t$ of $G$.

\begin{algorithm}
	\caption{}\label{v-d}
	\begin{algorithmic}[1]
		\State $H\gets G$
		\State $j \gets 0$
		\While{ $j\leq m$} 
		\State delete all nodes with degree smaller than $i-jt$ in $H$
		\State $H \gets \text{the resulting graph}$
		\State $j \gets j+1$
		
		\EndWhile
		\Return $H$
	\end{algorithmic}
\end{algorithm}

Consider the maximal D-chain of order $t$ of $G$, $G_0\geq G_1 \geq \cdots \geq G_k \geq \cdots$: a vertex is contained in $G_k$ iff at least $k$ of its neighbors are contained in $G_{k+t}$, i.e.~referencing
the vertex degree within $G_{k+t}$, a predecessor in the chain. This referencing propagates down to
$G_i$, after which we reference vertex degrees within $G$ itself. Reversing this backtracking, the
following vertex-deletion algorithm constructs the sequence $(G_i,G_{i-t},\dots,G_{i-mt}=G_k)$ starting
from $G$ as follows: it first constructs $G_i$, by deleting any vertices having $G$-degree less than
$i$ and then constructs $G_{i-t}\leq G_i$ by deleting all vertices of $G_i$-degree less than $i-t$
(note $t<0$). This continues inductively until it arrives at $G_k\leq G_{k+t}$, $k=m(-t)+i$.
The formal proof is given below.

\begin{theorem}\label{thm:deletion}
Let $G_0\geq G_1 \geq \cdots \geq G_k \geq \cdots$ be the maximal D-chain of order $t$ of the graph $G$.
Then the graph $H$ produced by Algorithm~\ref{v-d} equals $G_k$.
\end{theorem}
\proof
For any $l>0$, let $H_l$ be the graph produced by Algorithm~\ref{v-d}, for $k=l$. Firstly, for any $1\leq i \leq -t$,
the chain $L_i\colon H_0=G\geq H_i \geq H_{i-t}\geq H_{i-2t}\geq \cdots $ satisfies by the construction of Algorithm~\ref{v-d}
that any vertex $v$ contained in $H_{i-jt}$ has at least $i-jt$ neighbors in $H_{i-(j-1)t}$. 
Secondly, we claim that $L_i$ is maximal: there does not exist a chain
$L'_i\colon  H'_0=G\geq H'_i \geq H'_{i-t}\geq H'_{i-2t}\geq \cdots $
satisfying the same degree properties such that $H'_{i-jt}> H_{i-jt}$ for some $j\geq 0$.
To see this, suppose this were the case and let $j'$ be the minimal such that $H'_{i-j't}>H_{i-j't}$.
Then, there is a vertex $u$ contained in $H'_{i-j't}$ but not $H_{i-j't}$. However, this is impossible, since,
by assumption, we have $H'_{i-(j'-1)t}=H_{i-(j'-1)t}$ and the fact that $u$ is deleted from $H_{i-(j'-1)t}$ is tantamount to
$u$ having degree less than $i-j't$ in $H_{i-(j'-1)t}$. 

Finally, we claim that combining these $-t$ chains, that is, arriving at, $H_0, H_1, H_2, \ldots$, gives us the maximal D-chain of order
$t$ of $G$. It only remains to check the nesting property of adjacent graphs in the combined chain.
Suppose this is not true. Then there exists a minimum $r$ such that $H_r<H_{r+1}$.
Assume $r=i-jt$ whence $r+1=(i+1)-jt$. Note that both chains $L_i$ and $L_{i+1}$ satisfy the degree property and $H_{i-j't}\geq
H_{i+1-j't}$ for $0\leq j' < j$. Thus, the resulting chain obtained by replacing $H_{i-jt}$ with $H_{i+1-jt}$ in $L_i$ still
satisfies the degree property, contradicting the maximality just proved. This concludes the proof of Theorem~\ref{thm:deletion}. \qed

\subsection{The D-spectrum via fixed points}

In this subsection we present a different approach of computing the D-spectrum, namely, as a fixed point of a
discrete dynamical system. To this end, let us briefly recapitulate some basic facts about such systems. A
discrete dynamical system over a network involves the following
ingredients~\cite{kauf,vonn,wolf,rei1,rei2,rei5}:
a network, a local function associated with each node of the network that specifies how the
state of the node evolves and an update schedule that reflects when each individual node updates its state.
Von Neumann's cellular automata (CA) are
such systems. Given a network and local functions, the system dynamics is
concerned with how the system state varies in time. Various classes of dynamical systems have been
studied, e.g.~linear, sequential systems~\cite{linear1,linear2}, monotone systems~\cite{mon-neural1,mono1}, and
threshold systems~\cite{gole}. In particular, in Chen and Reidys~\cite{linear2}, a method of computing the M\"{o}bius function of a partially ordered set via implementing a discrete dynamical system was discovered. In the following, we shall see another instance of such usage of discrete dynamical systems for providing computation for problems arising from different fields.

Let $G=(V,E)$ be a network with vertex set $V=\{1,2,\ldots, n\}$ and edges in the set $E$. Suppose
each node, $i$, has states contained in a finite set $P$. We associate a function $f_{i}$, that specifies
how the vertex $i$ updates its state, $x_{i}$.
The update entails considering the states of the neighbors of $i$ and $i$ itself as arguments of $f_{i}$,
whence we call $f_{i}$ the local function at $i$.
An infinite sequence $W=W_1W_2\cdots$, where $W_i\subseteq V$, is called a fair update schedule, if
for any $k\geq 1$, and any $1\leq i \leq n$, there exists $l>k$ such that $i\in W_l$. The system dynamics
is being generated if nodes update their states using their respective local functions, following the order
specified by a fair update schedule $W$. That is, suppose the initial system state at time $t=0$ is
$x^{(0)}$. For $j>0$, the system state $x^{(j)}$ at time $t=j$ is obtained by the nodes contained in
$W_j$ updating their states by means of their local functions taking as arguments the states of their
respective neighbors in $x^{(j-1)}$. The states of the nodes not in $W_j$ remain unchanged.
We denote this dynamical system by $[G,f,W]$, and we write $[G,f,W]^{(j)}(x)$ for the system state at time
$t=j$ assuming the system has initial state $x$ at time $t=0$. For a given dynamical system $[G,f,W]$, a
system state $x$ is said to reach a fixed point (or stable state) $z$ if there exists $k\geq 0$ such that
for any $j>k$, we have $[G,f,W]^{(j)}(x)=[G,f,W]^{(k)}(x)=z$.

Suppose there is a linear order `$\leq $' on the set $P$. Let $P^q=\{(x_1,x_2,\dots, x_q)\mid x_j \in P,\, 1
\leq j \leq q\}$. We extend the linear order on $P$ to a partial order on $P^q$ as follows:
$(x_1,x_2,\dots , x_q)\leq (y_1,y_2,\dots, y_q)$ iff for all $1\leq j \leq q$, $x_j\leq y_j$ in $P$. A function
$g\colon P^q\rightarrow P$ is called monotone if for any $x\leq y$ in $P^q$, $g(x)\leq g(y)$ in $P$.
A local function $f_{i}\colon (x_i, x_{k_1},x_{k_2},\dots, x_{k_i})\mapsto x'_i$ is called contractive if
for any argument $(x_i, x_{k_1},x_{k_2},\dots, x_{k_i})\in P^{k_i+1}$, $x'_i \leq  x_i$.
For example, the Boolean functions `AND' and `OR' on $P^q=\{0,1\}^q$ are monotone, under both assumptions that $0<1$ and that $1<0$.
It is also easy to check that if $f_{v_i}$ is the Boolean function `AND', then it is contractive under the assumption $0<1$; while if $f_{v_i}$ is the Boolean function `OR', then it is contractive under the assumption $1<0$.
A dynamical system in which local functions are monotone and contractive is called a monotone-contractive
(MC) system. A key property of MC systems is the following:

\begin{theorem}\label{main-thm-1}
 For any two fair update schedules $W$ and $W'$, a system state $x\in P^n$ reaches the same fixed point $z\in P^n$ under the two MC systems $[G,f,W]$ and $[G,f,W']$. In addition, any state $y$ such that $z\leq y\leq x$ will reach the fixed point $z$.
\end{theorem}

The proof of Theorem~\ref{main-thm-1} is given in the Supplementary Information Section D. In addition, since from Theorem~\ref{main-thm-1} the exact update schedule does not matter, we will not explicitly specify the update schedule if not necessary in the following.  

Now we consider some concrete MC systems (w.r.t.~the part we are interested) on $G$ that have a close connection to the maximal D-chains of $G$.
Suppose each node has $[n]=\{1,2,\ldots, n\}$ as the set of states. Suppose the local function $f_v$ at node $v$ returns the maximum $k$ such that there are at least $k$ neighbors of $v$ having states at least
$k+t$. We call this system the $[t]$-system on $G$. Then the maximal D-chain of order $t$ and node specific
spectrum of $G$ can be computed based on Theorem~\ref{thm:main2}.

\begin{theorem}\label{thm:main2}
For the $[t]$-system on $G$, the state
$x=\big(\Deg(1),\Deg(2),\dots, \Deg(n)\big)$ reaches the stable state
$C^t=(C_{1,t},C_{2,t},\ldots, C_{n,t})$, where
	\begin{itemize}
		\item[i.] $C_{i,t}=0$ for any $1\leq i \leq n$, if $t>0$;
		\item[ii.] $C_{i,t}= C_t(i)$ for any $1\leq i \leq n$, if $t\leq 0$. In particular, if
                  $t\leq -\Delta(G)$, $C_{i,t}=\Deg(i)$ for any $1\leq i \leq n$.
	\end{itemize}
In addition, the state $C^t$ reaches the stable state $C^{t+1}$ in the $[t+1]$-system on $G$.
\end{theorem}
\proof 
From the theory of MC systems developed in the Supplementary Information (see the Supplementary Information Proposition~D.$1$, Proposition D.$5$), it follows that $x$ reaches a fixed point.
Suppose first, $t>0$ and suppose there exists some $v$ such that $C_{v,t}>0$. By definition, there is at least one neighbor
$u$ of $v$ such that $C_{u,t}\geq C_{v,t}+t>0$ holds. Iterating this argument, the node $u$ has a neighbor with $C_t$-value
at least $C_{u,t}+t$, effectively implying the existence of vertices with unbounded degrees, which is, given the fact that
the nework is finite, impossible.

Secondly, let $t\leq 0$.  We shall first prove that the state $z=\big( C_t(v_1), \ldots, C_t(v_n)\big)$ is a fixed point of
the $[t]$-system, i.e.~applying the local function $f_{v}$ (with parameter $t$) to $z$ will return $C_t(v)$ for
any vertex $v$. Let $L\colon G_0\geq G_1\geq \cdots $ be the maximal D-chain of order $t$ of $G$.
Then, for any vertex $v$, by definition of $C_t(v)=i$, $v$ belongs to $G_i$ but not to $G_{i+1}$. This implies the following:
\begin{itemize}
	\item[ (i)] there are at least $i$ neighbors of $v$ are contained in $G_{j}$ ($j=\max\{0, i+t\}$). Note that for any $u$
	among these, we have $C_t(u)\geq i+t$ by definition. Thus, among the neighbors of $v$, there are at least $i$ having
	values at least $i+t$ in $z$. This implies $f_v(z)\geq i=C_t(v)$;
	\item[(ii)] there cannot be at least $i+1$ neighbors of $v$, that are contained in $G_{j'}$ ($j=\max\{0, i+1+t\}$),
	as otherwise, the chain $G_0\geq \cdots \geq G_i\geq G_{i+1}\bigcup \{v\}\geq G_{i+2}\geq \cdots$ gives a D-chain of order $t$,
	which contradicts the maximality of $L$.
	Hence, among the neighbors of $v$, there can not be that at least $i+1$ of them with values at least $i+1+t$ in $z$, whence
	$f_v(z)< i+1$ holds.
\end{itemize}
(i) and (ii) establish $f_v(z)=i=C_t(v)$ and $z$ is a fixed point.

We proceed with the proof by observing that, in case of $z=x$, we are done. By construction, we otherwise have $z<x$.
Let $y$ be the fixed point reached by $x$. In case of $y<z$ or $y$ being
incomparable to $z$, Proposition~D.$2$ (in the Supplementary Information) guarantees that $z$ cannot be a fixed point, which is a contradiction.
Otherwise we have $z < y \leq x$. In this case there exists a coordinate, which we shall index by $v$, that
satisfies $y_v> C_t(v)$.
Consider the sequence of
subgraphs induced by the sequence of sets of vertices $S^0, S^1,\ldots$ which are inductively defined as follows:
(i) $S^0=\{v\}$;
(ii) for $r>0$, $S^r=\{u\mid  y_u\geq y_w+t, w \in S^{r-1}, \mbox{and $u=w$ or $u$ is a neighbor of $w$}\}$.
Clearly, by construction we have $S^{r-1}\subseteq S^r$, and by induction $y_u\geq y_v+rt$ for $u\in S^r$.
If $y$ is a fixed point, then we have: (a) there are at least $y_v$ neighbors of $v$ with values at least $y_v+t$
in $y$. These neighbors must be contained in $S^1$; (b) for $r\geq 0$, any $w\in S^r$, there are at least $y_w+t$
neighbors contained in $S^{r+1}$.

By abuse of notation, we will denote the subgraph induced by the set $S^r$-vertices as $S^r$. From (a), (b) and the fact
that $y_u\geq y_v+rt$ for $u\in S^r$, we can conclude that for $r\geq 0$, any vertex in $S^r$ has at least $y_v+rt$ neighbors in
$S^{r+1}$. Then, the chain 
\begin{multline*}
\cdots \geq G_{y_v+2t}\bigcup S^2 \geq G_{y_v+2t+1}\bigcup S^1 \geq \cdots\geq G_{y_v+t}\bigcup S^1 \\
\geq G_{y_v+t+1}\bigcup S^0 \geq \cdots \geq G_{y_v-1}\bigcup S^0\geq G_{y_v}\bigcup S^0\geq G_{y_v+1}\geq \cdots
\end{multline*}
is a D-chain of order $t$ of $G$, implying $C_t(v)\geq y_v$, which is a contradiction. Thus $x$ cannot reach a fixed point
$y$ such that $z<y\le x$, whence $x$ is reaching the fixed point $z$ as claimed. 

In particular, if $t\leq -\Delta(G)$, we have $C_t(v)=\Deg(v)$, whence $x$ is a fixed point.
Finally, Proposition~D.$5$ (in the Supplementary Information), implies $C^{t+1}\leq C^t \leq x$. We have just proved that the degree sequence,
$x$, converges to $C^{t+1}$ in a $[t+1]$-system on $G$. Theorem~\ref{main-thm-1} in turn implies that since $C^{t+1}\le C^t$, $C^t$ also converges to $C^{t+1}$ in the $[t+1]$-system.\qed

\emph{Remark.}
The fact that the state $C^t$ converges to the stable state $C^{t+1}$ for the $[t+1]$-system on $G$, as claimed in Theorem~\ref{thm:main2}, guarantees essentially the same complexity for computing the D-spectra of all nodes as computing core numbers alone,
that is, it is not necessary to start with the degree sequence every time.

\begin{equation*}
\xymatrix{x\ar[d]_{-\Delta(G)}\ar@/^4pc/[ddddd]^0\\
	C^{-\Delta(G)}\ar@{.>}[d]\\	
	C^t\ar[d]_{t+1}\\
	C^{t+1}\ar[d]_{t+2}\\
	C^{t+2}\ar@{.>}[d]\\
	C^0
}
\end{equation*}


\section{Application to predicting similarity}


 Here we shall be
interested in analyzing the connection between the D-spectra and the spreading power of nodes in the process such as disease outbreak or information spreading. Specifically, we will be using the SIR model to get the data on infection rates charaterizing the spreading power of nodes.  
To begin with, let us briefly review the SIR model and our simulatin setup.
The SIR process is a stochastic model for studying the spread of disease within a population. It works as follows: a
population is modeled as a network, where each node represents an individual, while links (edges) between nodes represent their
interaction relation. Each node can be in either of three states: susceptible (S), infected (I), and recovered (R).
During the process, at each step, each infected node may infect each of its susceptible neighbors with a certain probability. At
the subsequent step the infected node may become recovered with another probability. Once a node is in the state R, it will never
infect other nodes and never become infected again. The process stops when there are no nodes in the state I.

Our SIR simulations are designed as follows: for each respective network, we initialize the process with exactly one node,
the infected source, in the state I. We shall assume that the probability of an infected node becoming recovered in the next
time step equals $1$ and we assume one fixed transmission probability for all nodes throughout the simulation. For each infected
source, we run $1000$ simulations and for each of these we compute the ratio between the number of recovered nodes and the total
number of nodes in the network. We refer to the average of these ratios as the infection rate of the node.

We execute this for each node in the network, for the nine transmission probabilities $h\cdot \beta$, where
$h\in \{0.1,0.5,1, 1.5, 2,4,6,8,10\}$ and $\beta$ being the epidemic threshold value of the network.
The epidemic threshold value can be computed as $\beta=\frac{<k>} {<k^2>-<k>}$~\cite{spread1,spread2}, where $<k>$ denotes the
average node degree of the network and $<k^2>$ is the average of the squares of the degrees. Accordingly, we obtain for each node nine distinct infection rates.

Kitsak et al.~\cite{np-core} shows that with respect to infection rates, nodes that are contained in
the same core are generally more isotropic, than nodes having the same degree. In other words, in
order to identify nodes of similar spreading power, core numbers are more suited than vertex degrees. 
In the following, we compare D-spectrum and core number.

Recall that the D-spectrum of a node is a vector whose coordinates are the ranks of the node in the respective D-chains of order $t$.
As a result, the Euclidian distance between D-spectra is a natural criterion for categorizing (possibly) similar nodes. In order to compare such
a categorization to the one derived by restricting to core numbers, i.e.~nodes having the same core number being considered
similar, we study five real networks and for each network proceed as follows. We partition all nodes by means of the Euclidean distance between their D-spectra into an
{\it a priori} specified number of clusters called D-blocks. (The clusters are derived by calling the standard function
Findcluster in Mathematica 10.0.) We also separately group the nodes by their core numbers into clusters that we refer to as C-blocks.

By construction, a D-block may contain nodes from multiple C-blocks and vice versa. A pairwise intersection of C- and
D-blocks is called an I-cell. The I-cells provide the underlying grid of Figure~$2$, (a), (c), (e), where the
color of the corresponding I-cell represents the average of the infection rates of the nodes contained in that cell (at
the same fixed transmission probability). In Figure~$2$, (b), (d), (f), in order to quantify
the comparison of the partitions of the nodes by D- or C-blocks, first we compute for each
C-block the dispersion (variance-to-mean ratio) of infection rates of nodes within said block. Secondly, we compute the dispersion of the
infection rates of nodes within each of the I-cells that refine said C-block. These data, provide a local to global picture
that quantifies the refinement that D-spectra provide over conventional cores.

The studied five networks are: Email~\cite{email},
USAir~\cite{usair}, Jazz~\cite{jazz}, PB~\cite{pb}, and Router~\cite{router}, described in detail
in the Supplementary Information. For the analysis of the three networks of Figure~$2$, the transmission probabilities
were set to $1.5  \beta$, where $\beta$ is their respective epidemic threshold value. 

\begin{figure}[H]
\label{F:22}
	\setlength{\tabcolsep}{5pt}
	\centering
	\begin{tabular}{p{0.55\textwidth}p{0.35\textwidth}}
		\begin{subfigure}{0.55\textwidth}
				\caption{Email}
			\centering
			\includegraphics[width=0.88\linewidth]{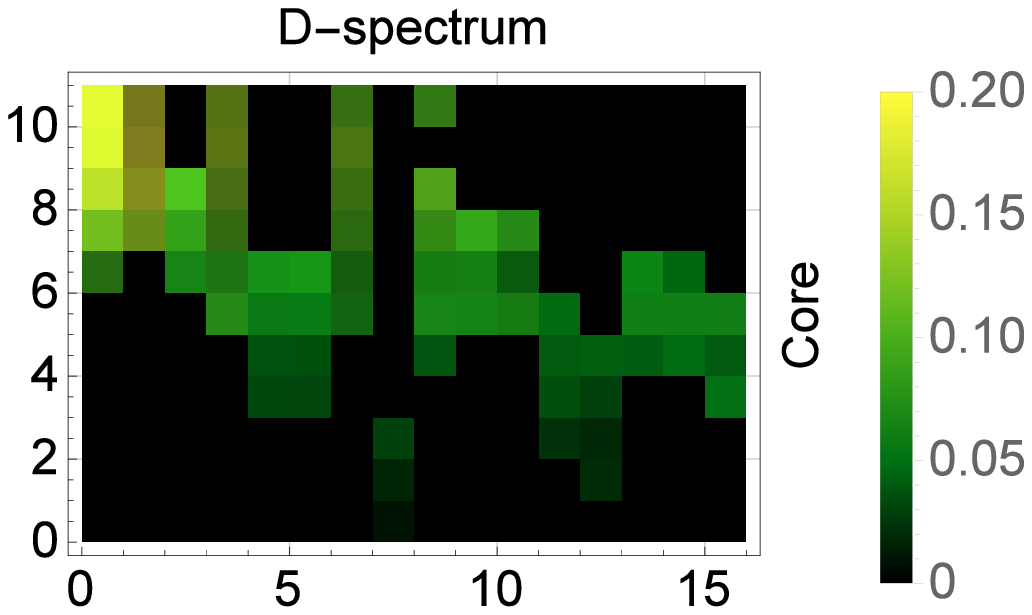}
		\end{subfigure}%
		&	\begin{subfigure}{0.35\textwidth}
				\caption{Email-dispersion}
			\centering
			\includegraphics[width=1.0\linewidth]{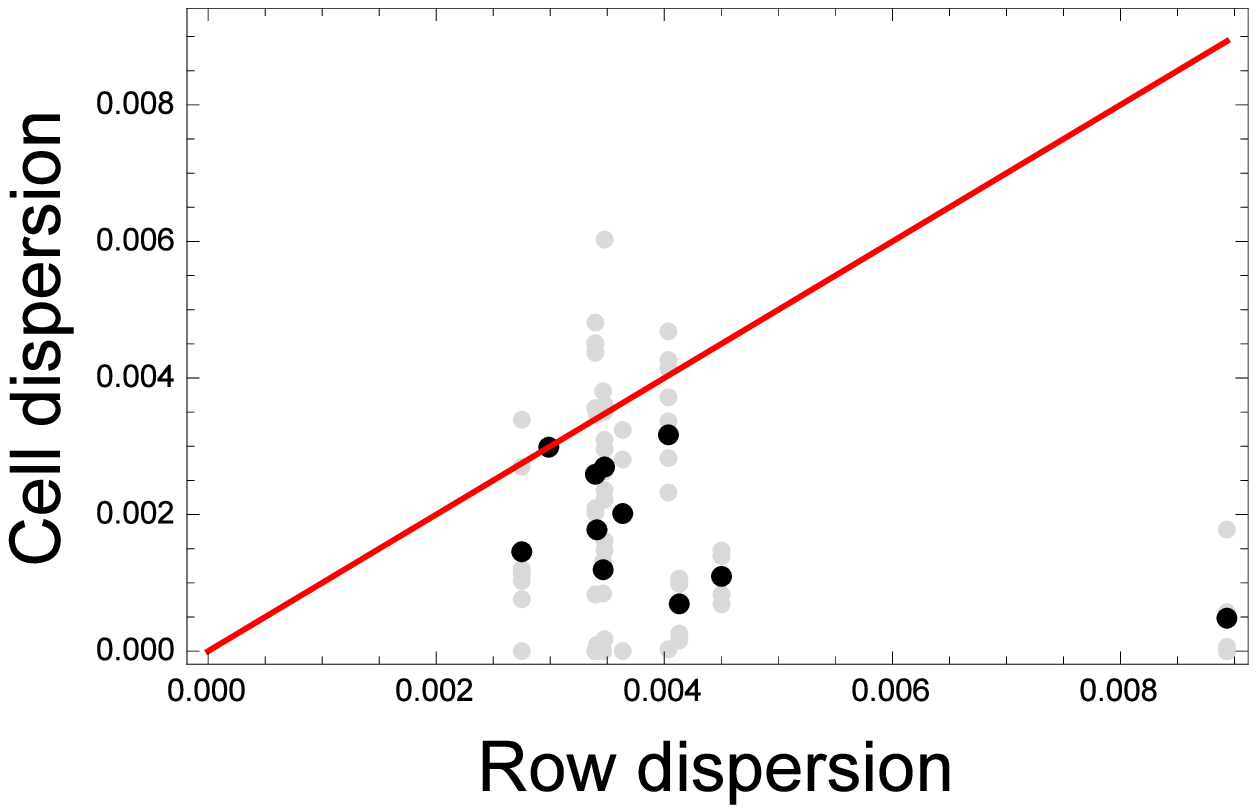}
		\end{subfigure}\\
		\begin{subfigure}{0.55\textwidth}
				\caption{Jazz}
			\centering
			\includegraphics[width=0.88\linewidth]{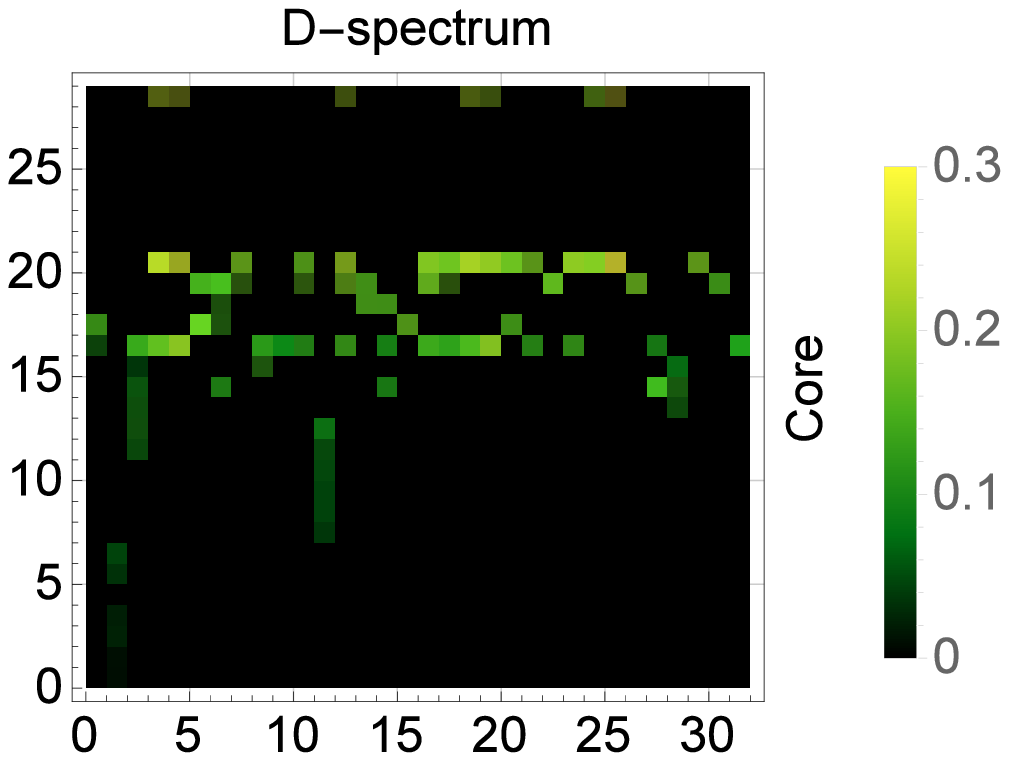}
		
		\end{subfigure}%
		& 	\begin{subfigure}{0.35\textwidth}
				\caption{Jazz-dispersion}
			\centering
			\includegraphics[width=1.0\linewidth]{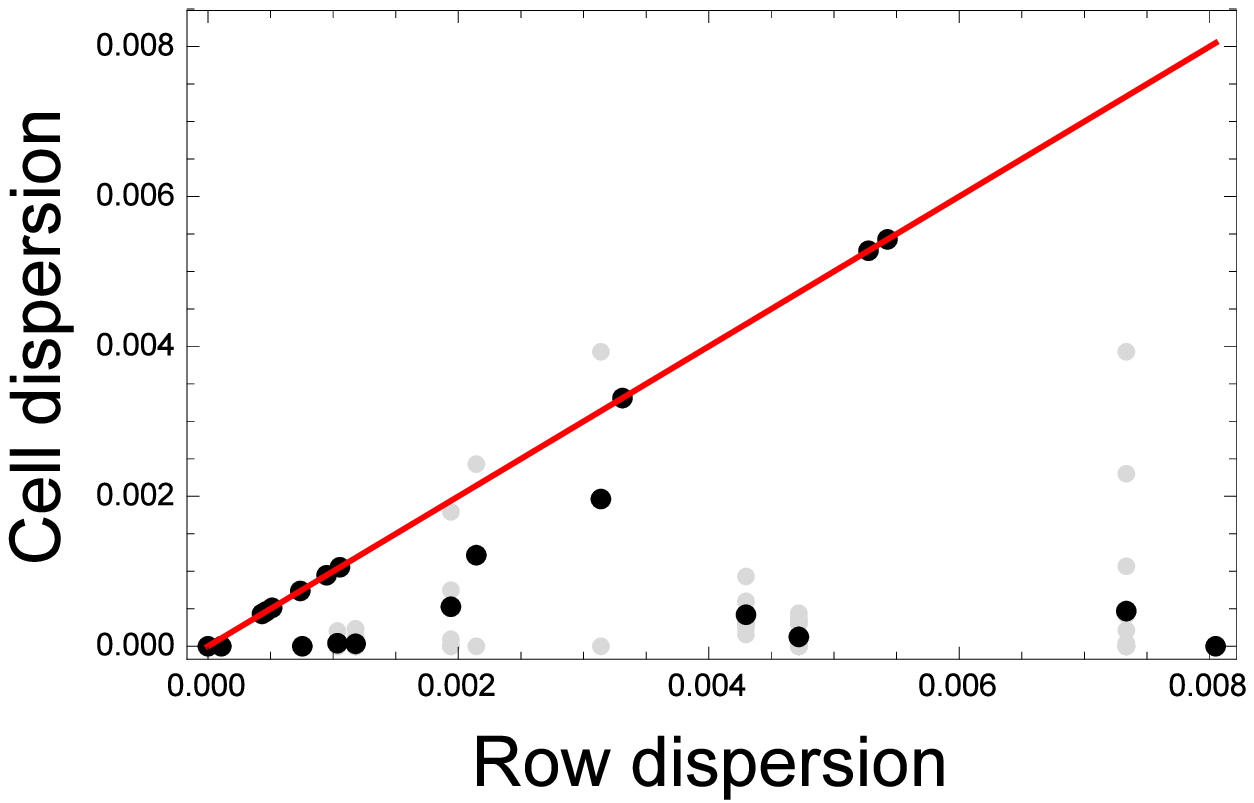}
		
		\end{subfigure}\\
		\begin{subfigure}{0.55\textwidth}
				\caption{PB}
			\centering
			\includegraphics[width=0.88\linewidth]{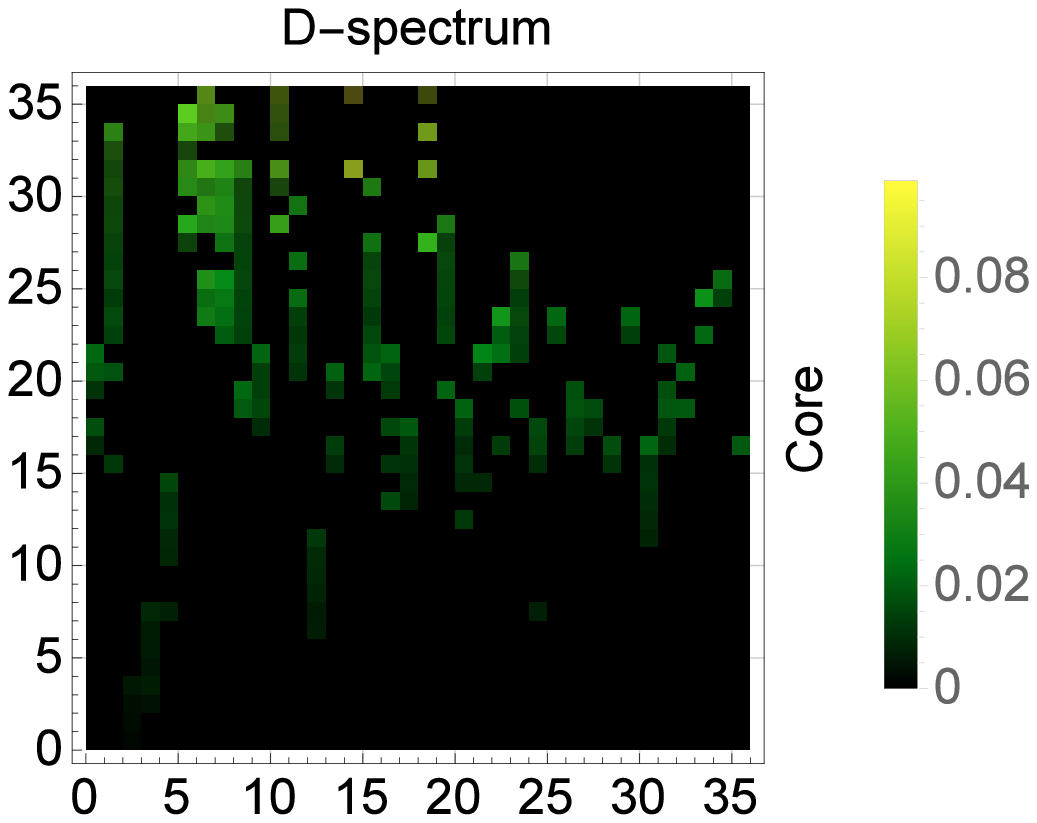}
		
		\end{subfigure}
		&	\begin{subfigure}{0.35\textwidth}
				\caption{PB-dispersion}
			\centering
			\includegraphics[width=1.0\linewidth]{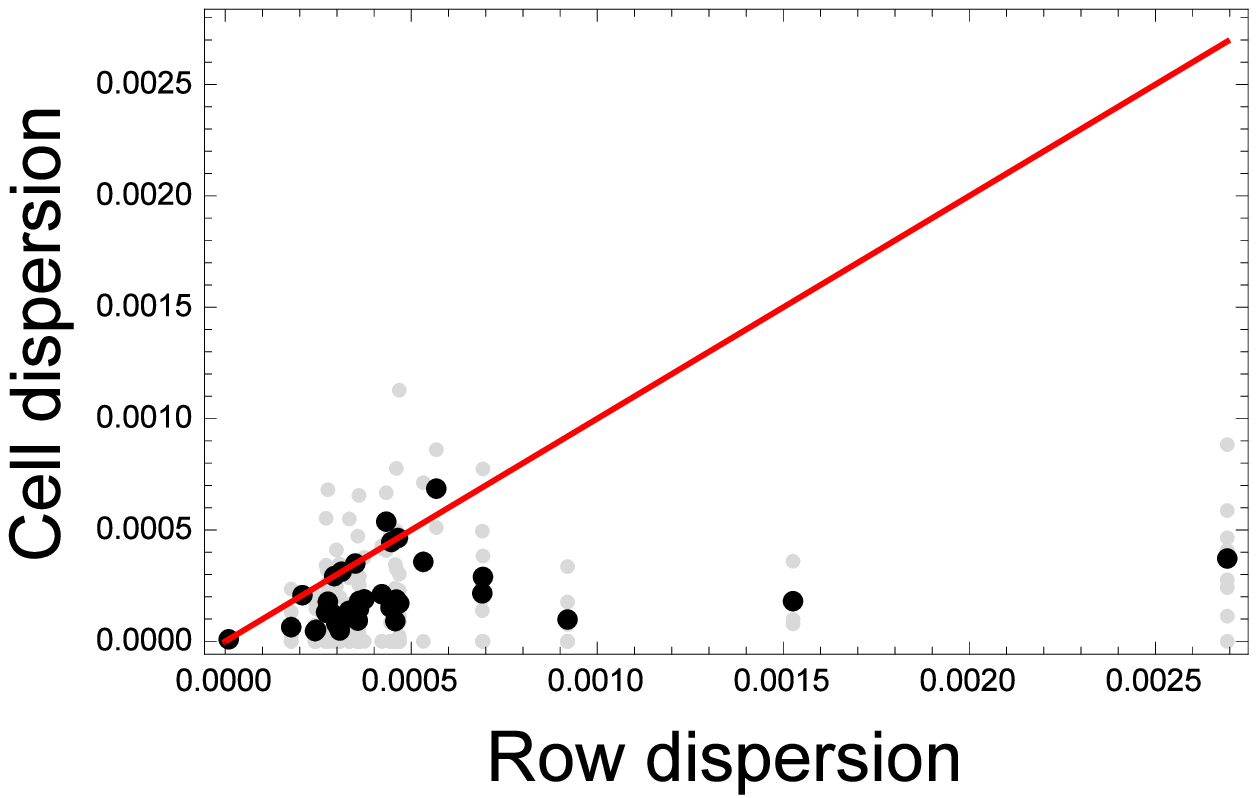}
		
		\end{subfigure}%
	\end{tabular}
	\caption{{\bf D-spectra versus cores.}
          {\normalfont In (a),(c),(e), rows and columns represent C- and D-blocks, respectively, while the color of each
            I-cell corresponds to the average of the infection rates of nodes contained in the I-cell at 
            the transmission probability $1.5\beta$. For all networks, the average infection rates in the same D-block (column) are more isotropic than those
            of the same C-block (row).
            Figures (b),(d),(f) provide a quantification of the local versus global dispersions:
            the x-coordinate represents the global dispersion of a given C-block. The $y$-coordinates (gray), for a fixed $x$,
            correspond to the dispersions for the I-cells that refine the C-block corresponding to $x$.
            The $y$-coordinate (black) represents the average of the I-cell dispersions for each $x$ fixed.}}
\end{figure}

In the Supplementary Information, we extend the analysis of these networks incorporating the following additional two transmission probabilities: $1\beta$ and $2\beta$.  
The analysis of the additional probabilities shows the robustness of the observation from Figure~$2$ that the infection rates of nodes
having the same core number, are generally more heterogeneous than that of nodes in the same D-spectrum block. See Figure~$2$ and Figures S1--S13 in the Supplementary Information. Accordingly, node partitions obtained via D-spectra provide a meaningful enhancement over
categorizations obtained using conventional cores.

We next qualify the correlation between the spreading power of nodes observed in the SIR process and the D-spectra of
nodes. That is, we ask to what extent do nodes, categorized via D-spectra, exhibit isotropic spreading power in the SIR
process. To this end, we firstly cluster the nodes according to their spreading power (i.e.~the sequences of infection rates at the nine transmission probabilities) and secondly we cluster them w.r.t.~their
D-spectra. 
Then we inspect the mutual intersection of these clusters from the two approaches.
We find that these two partitions are highly correlated in most cases and the correlation is robust w.r.t.~different specified number of clusters: see Figure~$3$ as well as Figures S14--S18 in the Supplementary Information.

\begin{figure}[!htb]
	\setlength{\tabcolsep}{5pt}
	\centering
	\begin{tabular}{p{0.28\textwidth}p{0.33\textwidth}p{0.38\textwidth}}
		\begin{subfigure}{0.28\textwidth}
				\caption{Email}
			\centering
			\includegraphics[width=1.0\linewidth]{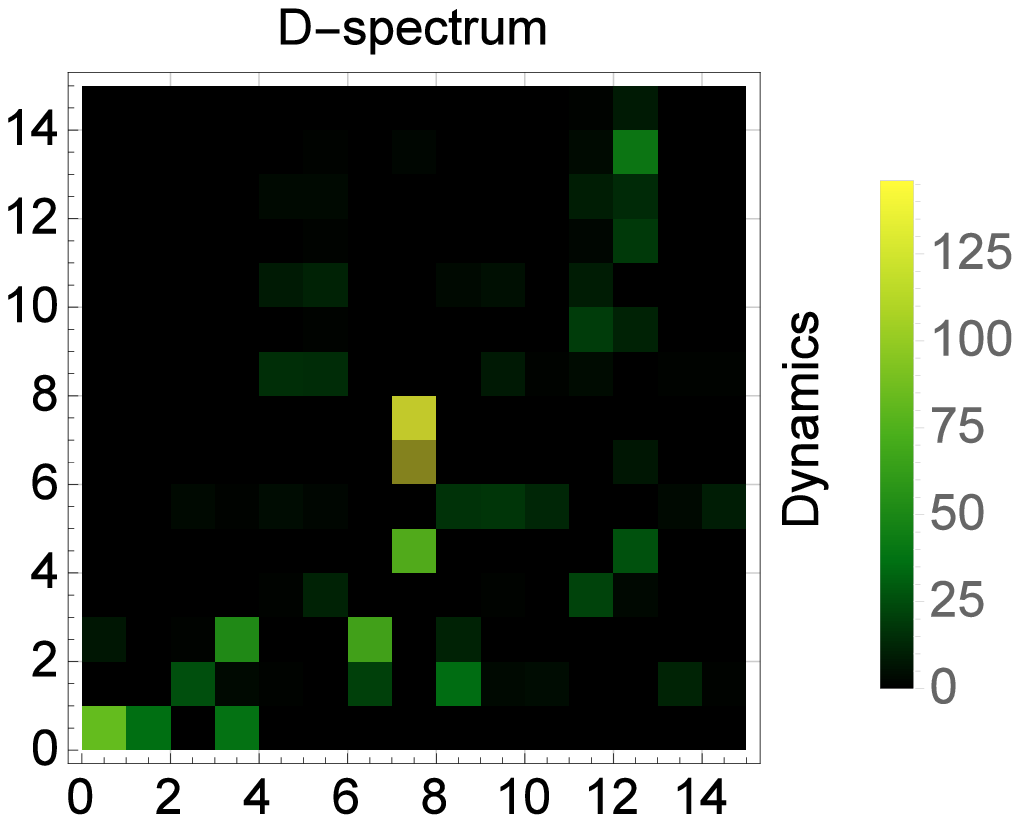}
		
		\end{subfigure}%
		&	\begin{subfigure}{0.33\textwidth}
				\caption{Jazz}
			\centering
			\includegraphics[width=1.0\linewidth]{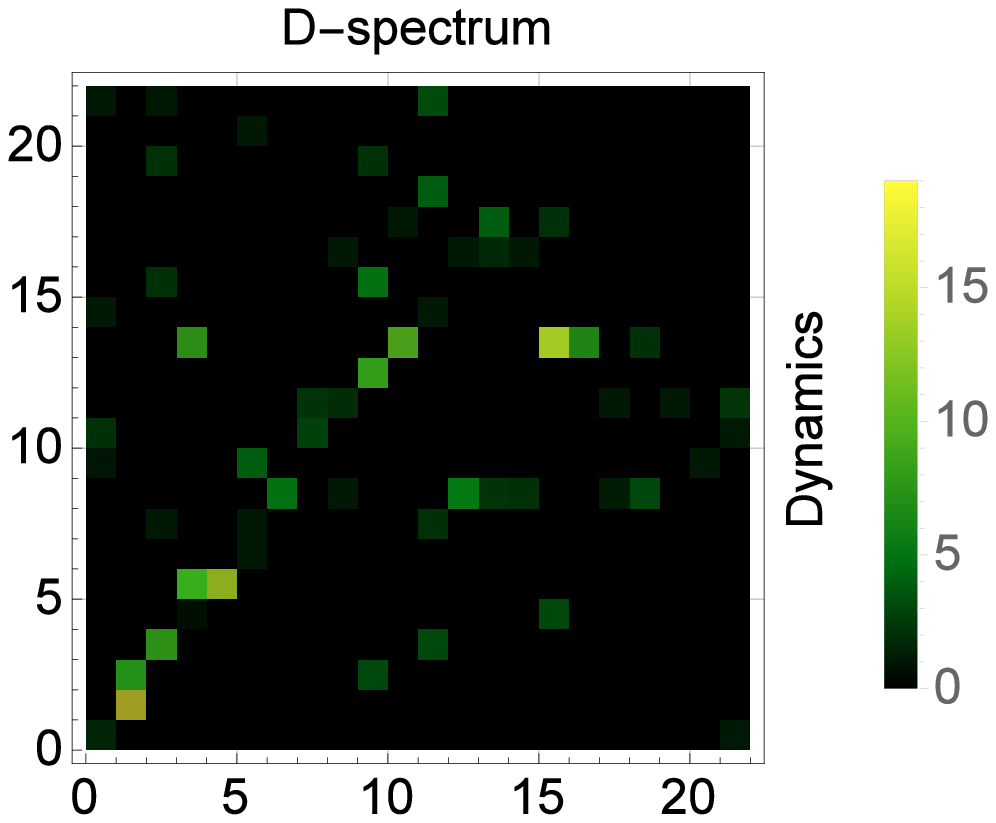}
		
		\end{subfigure}%
		&	\begin{subfigure}{0.38\textwidth}
				\caption{PB}
			\centering
			\includegraphics[width=1.0\linewidth]{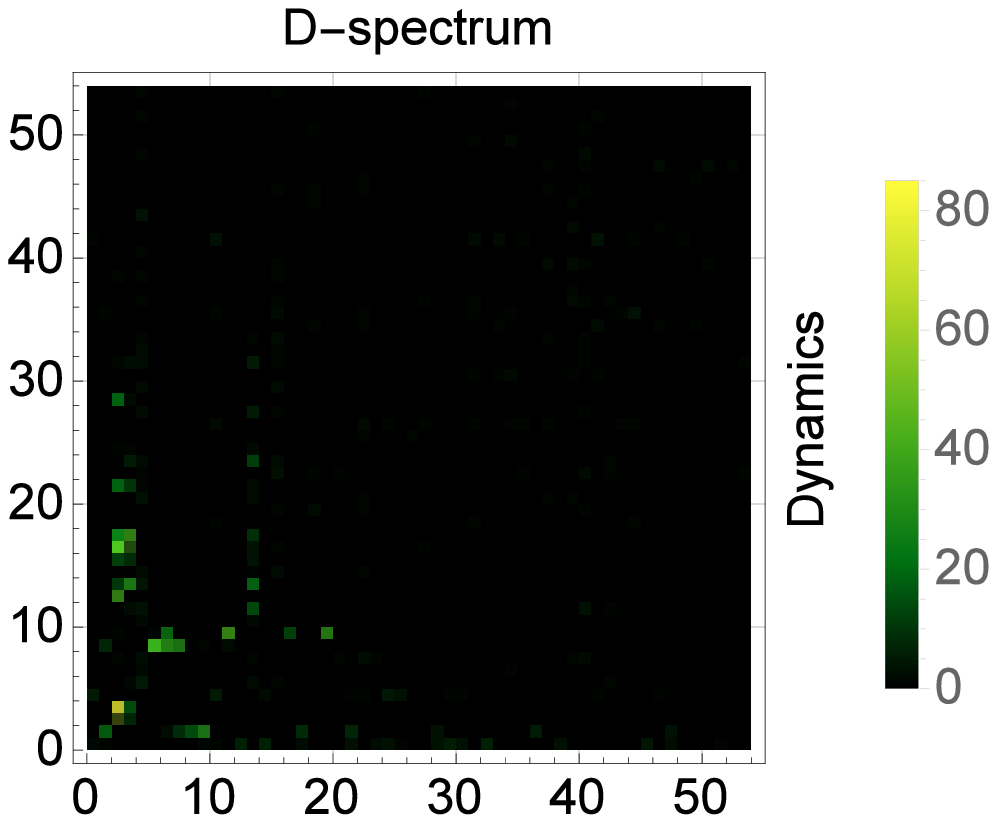}
		
		\end{subfigure}%
	\end{tabular}
	\caption{{\bf Partitions induced by D-spectra and spreading power are highly correlated.}
          {\normalfont In (a),~(b),~(c), the color of any cell represents the size of the intersection between
            blocks of the partitions induced by spreading power and D-spectra, respectively. For all networks, we observe
            that there are only a few distinct cells for any row or column that contain almost all vertices. The extent of this
            concentration reflects how well D-spectra capture the spreading power of vertices.}}
\end{figure}

\section{Discussion}
 
In this paper we develop a theoretical framework based on D-chains of certain levels of a given network. We establish uniqueness of the maximal D-chains
and show how to compute all relevant data via a parametric deletion algorithm as well as by computing the fixed point of a particular
graph dynamical system, to which we refer as an MC system. The framework itself is not restricted to graphs, and it can be extended
to hyper-graphs, weighted networks, and $k$-truss decompositions~\cite{truss1,truss2}.
Varying our adaptive parameter $t$ we endow each node of the network with a vector of data called its D-spectrum. Compared to the
degree or the core number of a node, both of which appear as particular coordinates in this associated vector, the D-spectrum of
the node is a high dimensional object. As such it provides more information about the node and offers more flexibility when it comes
to its use as a classifier. This is illustrated by the fact that there are no natural criteria for further classifying nodes with the
same degree or core number. The resolution of classification of the nodes of a network is thus {\it a priori} determined by the number
of different degrees or cores respectively. In contrast, the D-spectra are much more refined and thus facilitate the separation of any two nodes, at higher resolution.

We then apply D-spectra to categorize nodes in a network. D-spectra can naturally be partitioned using the Euclidean distance. 
We then show that D-spectra are a good predictor of node spreading power within the context of SIR dynamics on the network, as well
as being by construction a good measure of structural similarity. The latter gives rise to the question to what extent two graphs,
exhibiting the same D-spectrum are similar.
 
The framework presented here is far from being fully explored. This holds for theory as well as applications. For instance, are two
networks having the same D-spectrum (instead of a single ranking such as degree or core number) isomorphic? As for applications we restricted our analysis to SIR processes, which are well
studied and in which the relevancy of degree and core-number have been established. Since the theme here is to connect network
structure with dynamics, it is worth studying other processes such as, voting, information dissimination and processes whose
transmissions are mediated by other threshold functions. It can be speculated that, depending on the process, for certain $t$,
D-chains of order $t$ are more relevant classifiers than others.

\section*{Acknowledgements}
We thank Stephen Eubank and Henning Mortveit for valuable discussions.
We also thank Linyuan L\"{u} and Qian-Ming Zhang for providing some data related to the used networks.

\section*{Author information}
\subsection*{Affiliations}
Biocomplexity Institute and Initiative, University of Virginia, Charlottesville, Virginia 22908, USA\\
Ricky X. F. Chen \& Christian M. Reidys

\noindent Department of Mathematics, Virginia Tech, Blacksburg, Virginia 24601, USA\\
Andrei C. Bura

\subsection*{Contributions}

R.X.F.C.~and C.M.R.~planned and performed this research. R.X.F.C.~partly implemented the simulation. A.C.B.~implemented the simulation and performed the research. All authors discussed the results, wrote the paper and reviewed the manuscript.

\subsection*{Competing interests}
The authors declare no competing financial interests.

%
%
%
%
%
%
%
\subsection*{Corresponding authors}
Correspondence to Ricky X. F. Chen (chen.ricky1982@gmail.com) or Christian M. Reidys (duck@santafe.edu).


\begin{thebibliography}{99}
	
	\bibitem{scale} Barabasi, A.-L. \& Albert, R. Emergence of scaling in random networks. {\it Science}
	{\bf 286,} 509–512 (1999).
	
		\bibitem{scale2} Albert, R., Jeong, H. \& Barabasi, A.-L. Error and attack tolerance of complex
	networks. {\it Nature} {\bf 406,} 378–382 (2000).
	
	\bibitem{path} Watts, D. J. \& Strogatz, S. H. Collective dynamics of `small-world' networks.
	{\it Nature} {\bf 393,} 440–442 (1998).
	
	\bibitem{path2} Freeman, L. C. A set of measures of centrality based on betweenness.
	{\it Sociometry} {\bf 40,} 35–41 (1977).
	
	\bibitem{kcore1} Dorogovtsev, S. N., Goltsev, A. V. \& Mendes, J. F. F. K-core organization of
	complex networks. {\it Phys. Rev. Lett.} {\bf 96,} 040601 (2006).
	
		\bibitem{kcore2} Seidman, S. B. Network structure and minimum degree. {\it Social Networks} {\bf 5,} 269–287 (1983).
	
	\bibitem{kcore3} Carmi, S., Havlin, S., Kirkpatrick, S., Shavitt, Y. \& Shir, E. A model of Internet topology using k-shell decomposition. {\it Proc. Natl Acad. Sci. USA} {\bf 104,} 11150–11154 (2007).
	
	\bibitem{motif} Alon, U. Network motifs: theory and experimental approaches. {\it Nat. Rev. Genet.}
	{\bf 8,} 450–461 (2007).
	
	\bibitem{commu1}Newman, M. E. J. \& Girvan, M. Finding and evaluating community structure in
	networks. {\it Phys. Rev. E} {\bf 69,} 026113 (2004). 
	\bibitem{commu2} Newman, M. E. J. Finding community structure in networks using the
	eigenvectors of matrices. {\it Phys. Rev. E} {\bf 74,} 036104 (2006).
	
	\bibitem{spread1} Castellano, C. \& Pastor-Satorras, R. Thresholds for epidemic spreading in
	networks. {\it Phys. Rev. Lett.} {\bf 105,} 218701 (2010). 
	
	\bibitem{spread2} Newman, M. E. J. Spread of epidemic disease on networks. {\it Phys. Rev. E} {\bf 66,}
	016128 (2002).
	
	\bibitem{spread3} Pastor-Satorras, R. \& Vespignani, A. Epidemic spreading in scale-free networks.
	{\it Phys. Rev. Lett.} {\bf 86,} 3200 (2001).
	
	\bibitem{spread4} Keeling, M. J. \& Rohani, P. Modeling Infectious Diseases in Humans and
	Animals (Princeton Univ. Press, 2008).
	
	
	
	\bibitem{spread5} Hethcote, H. W. The mathematics of infectious diseases. {\it SIAM Rev.} {\bf 42,}
	599–653 (2000).
	
		\bibitem{spread6} Eubank S. et al. Modelling disease outbreaks in realistic urban social networks. {\it Nature} {\bf 429,} 180–184 (2004).
		
	\bibitem{spread7} Rogers, E. M. {\it Diffusion of Innovation} 4th edn (Free Press, 1995).
	
		\bibitem{import1} Pastor-Satorras, R. \& Vespignani, A. Immunization of complex networks. {\it Phys.
	Rev. E} {\bf 65,} 036104 (2002).
	
	\bibitem{import2} Wang, P., Lu, J. \& Yu, X. Identification of important nodes in directed
	biological networks: a network motif approach. {\it PLoS ONE} {\bf 9,} e106132 (2014).
	


	
	\bibitem{import3} Morone, F. \& Makse, H. A. Influence maximization in complex networks
	through optimal percolation. {\it Nature} {\bf 524,} 65–68 (2015).
	

	

	
	\bibitem{import5} Anderson, R. M., May, R. M. \& Anderson, B. {\it Infectious Diseases of Humans: Dynamics and Control} (Oxford Science Publications, 1992).
	
	
	
	\bibitem{import6} Aral, S. \& Walker, D. Identifying Influential and Susceptible Members of Social Networks. {\it Science} {\bf 337,} 337-341 (2012).
	

	
	
		\bibitem{np-core} Kitsak, M. et al. Identification of influential spreaders in complex networks.
	{\it Nat. Phys.} {\bf 6,} 888–893 (2010).
	
		\bibitem{kcore-b} Bollobás, B. {\it Graph Theory and Combinatorics: Proceedings of the Cambridge Combinatorial Conference in Honor of P. Erdös} Vol. 35 (Academic, 1984).
	

	
	\bibitem{nc-core} L\"{u}, L., Zhou, T., Zhang, Q.-M. \& Stanley, H.E. The h-index of a network node and its relation to degree and coreness. {\it Nat. Commun.} {\bf 7,} 10168 (2016).
	

	
		\bibitem{mont} Montresor, A., Pellegrini, F. D. \& Miorandi, D. Distributed k-core decomposition. {\it IEEE Trans. Parallel Distrib. Syst.} {\bf 24,} 288–300 (2013).
	
		\bibitem{vonn} Von Neumann, J. {\it Theory of Self-Reproducing Automata} (University of Illinois Press, Chicago, 1966).
	
		\bibitem{kauf} Kauffman, S. A. Metabolic stability and epigenesis in randomly constructed genetic nets. {\it J. Theor. Biol.} {\bf 22,} 437-467 (1969).
	
	
	\bibitem{wolf} Wolfram, S. {\it Cellular Automata and Complexity} (Addison-Wesley, New York, 1994).


	\bibitem{rei1} Barrett, C. L. \& Reidys, C. M. Elements of a theory of simulation I. {\it Appl. Math. Comput.} {\bf 98,} 241-259 (1999).
	\bibitem{rei2} Barrett, C. L., Mortveit, H. S. \& Reidys, C. M. Elements of a theory of simulation II: sequential dynamic systems.
	{\it Appl. Math. Comput.} {\bf 107,} 121-136 (2000).





\bibitem{rei5} Mortveit, H. S. \& Reidys, C. M. {\it An Introduction to Sequential Dynamic Systems} (Springer, 2008).

\bibitem{linear1} Elspas, B. The theory of autonomous linear sequential networks. {\it IRE Transactions on
Circuit Theory} {\bf 6,} 45-60 (1959).
	
	\bibitem{linear2} Chen, R. X. F. \& Reidys, C. M. Linear sequential dynamical systems, incidence algebras, and M\"{o}bius functions. {\it Linear Algebra Appl.} {\bf 553,} 270-291 (2018).
	
	\bibitem{mono1} Chen, R. X. F., Mortveit, H. S. \& Reidys, C. M. Dependence of update schedules of monotone sequential Boolean networks, submitted.
	
		\bibitem{mon-neural1} Daniels, H. \& Velikova, M. Monotone and Partially Monotone Neural Networks.
	{\it IEEE Trans. Neural Netw.} {\bf 21,} 906-917 (2010).
	
	
	

	

	
	

\bibitem{gole} Goles, E. {\it Comportement oscillatoire d'une famille d'automates cellulaires non uniformes} (Th\`{e}se
IMAG, Grenoble 1980).
	
	\bibitem{hirsch} Hirsch, J. E. An index to quantify an individual’s scientific research output.
{\it Proc. Natl Acad. Sci. USA} {\bf 102,} 16569-16572 (2005).
	



	



	\bibitem{email} Guimer\`{a}, R., Danon, L., Diaz-Guilera, A., Giralt, F. \& Arenas, A. Self-similar community structure in a network of human interactions. {\it Phys. Rev. E} {\bf 68,} 065103 (2003).

\bibitem{usair} Batageli, V. \& Mrvar, A. Pajek Datasets. Available at http://vlado.fmf.uni-lj.si/pub/networks/data/2007.

\bibitem{jazz} Gleiser, P. \& Danon, L. Community structure in Jazz. {\it Adv. Complex Syst.} {\bf 6,} 565 (2003).

\bibitem{pb} Adamic, L. A. \& Glance, N. in: Proceedings of the 3rd International Workshop on Link Discovery. 36–43 (ACM 2004).	

\bibitem{router} Spring, N., Mahajan, R., Wetherall, D. \& Anderson, T. Measuring ISP topologies with Rocketfuel. {\it IEEE/ACM Trans. Networking} {\bf 12,} 2–16 (2004).

\bibitem{email2} Email Dataset. Available at http://www-levich.engr.ccny.cuny.edu/webpage/hmakse\\
/software-and-data/

\bibitem{adult} The Internet Movie Database. Available at http://www.imdb.com. 







	

	


\bibitem{truss1} Cohen, J. Trusses: Cohesive subgraphs for social network
analysis. (2008).

\bibitem{truss2}  Wang, J. \& Cheng, J. Truss decomposition in massive networks. {\it Proceedings of the VLDB Endowment} {\bf 5,} 812–823 (2012).


\end{thebibliography}
\end{document}